\documentclass[letterpaper, 10 pt, conference]{ieeeconf}  

\IEEEoverridecommandlockouts                              

\usepackage{graphics} 
\usepackage{epsfig} 
\usepackage{amsmath} 
\usepackage{amssymb}  
\usepackage{caption}

\usepackage{verbatim}
\usepackage{graphicx} 
\usepackage{epsfig} 
\usepackage{epstopdf}
\usepackage{amsmath} 
\usepackage{amssymb}  
\usepackage{caption}
\usepackage{color}

\usepackage{mathtools}
\makeatletter
\newcommand{\figcaption}{\def\@captype{figure}\caption}
\newcommand{\tabcaption}{\def\@captype{table}\caption}

\makeatother
\title{\LARGE \bf At What Frequency Should the Kelly Bettor Bet?
}

 \author{\large Chung-Han Hsieh,$^{1}$ B. Ross Barmish,$^{2}$ and John A. Gubner$^{3}$
 	\thanks{\hskip -10pt ${}^1$Chung-Han Hsieh is a graduate student working towards to his Ph.D. degree in the Department of Electrical and Computer Engineering, University of Wisconsin, Madison, WI 53706. E-mail: hsieh23@wisc.edu.}
 	\thanks{\hskip -10pt ${}^2$B. Ross Barmish is a faculty member in  the Department of Electrical and Computer Engineering, University of Wisconsin, Madison, WI 53706. \mbox{E-mail}: barmish@engr.wisc.edu.}
 	\thanks{\hskip -10pt ${}^3$John A. Gubner is a faculty member in  the Department of Electrical and Computer Engineering, University of Wisconsin, Madison, WI 53706. \mbox{E-mail}: john.gubner@wisc.edu.}
}

\begin{document}

\maketitle
\thispagestyle{empty}
\pagestyle{empty}

\parindent = 0pt
\begin{abstract}
	We study the problem of optimizing the \mbox{\it betting frequency} in a dynamic game setting using Kelly's celebrated expected logarithmic growth criterion as the performance metric. The game is defined by a sequence of bets with independent and identically distributed returns~$X(k)$. The bettor selects the fraction of wealth~$K$  wagered at~$k = 0$ and waits~$n$ steps before updating the bet size. Between updates, the proceeds from the previous bets remain at risk in the spirit of ``buy and hold.'' Within this context, the main questions we consider are as follows: How does the optimal performance, we call it~$g_n^*$, change with~$n$? Does the high-frequency case,~$n = 1$, always lead to the best performance? What are the effects of accrued interest and transaction costs? First, we provide rather complete answers to these questions for the important special case when~$X(k) \in \{-1,1\}$ is a Bernoulli random variable with  probability~$p$ that~$X(k) = 1$. This  serves as an entry point for future research using a binomial lattice model for stock trading. The latter sections focus on more general probability distributions for~$X(k)$ and two conjectures. The first conjecture is simple to state: Absent transaction costs,~$g_n^*$ is non-increasing in~$n$. The second conjecture involves the technical condition which we call the \mbox{\it sufficient attractiveness inequality}. We first prove that satisfaction of this inequality is sufficient to guarantee that the low-frequency bettor using large~$n$ can match the performance of the high-frequency bettor using~$n=1$. Subsequently, we conjecture, and provide supporting evidence that this condition is also necessary.
\end{abstract}

\section{INTRODUCTION}
\label{Section: INTRODUCTION }
\vspace{-1mm}
The takeoff point for this research is Kelly's expected logarithmic growth criterion which was originally used as a performance metric in a variety of sequential betting problems; see~\cite{Kelly_1956} and further developments in~\cite{Latane_1959}-\cite{Thorp_1969}. This criterion is not only fundamental to gambling theory but also a starting point for a line of research on portfolio optimization in the stock market; e.g., see  \cite{Cover_Thomas_2012}-\cite{MacLean_Thorp_Ziemba_2011}. In addition, within this body of literature, several desirable properties resulting from use of the Kelly criterion have been established. Most notable of these properties are the asymptotic performance guarantees which are established in some of these papers. Finally, a sampling of more recent papers on these topics includes~\mbox{\cite{Kuhn_Luenberger_2010}-\cite{Hsieh_Barmish_Gubner_2016_CDC}.}

\vspace{3mm}
To be more specific, the classical Kelly Betting Problem is described by a sequence of bets with independent and identically distributed (i.i.d.) returns~$X(k)$ with known probability distribution. The bettor designates the fraction of wealth~$K$ being wagered seeking to maximize the expected logarithmic growth of the account value~$V(k)$ from $k=0$ until the terminal stage~$k=N$. At intermediate stage~$k$, if the account value is~$V(k)$, and since the bet size is~$KV(k)$, the gain or loss is given by~$KX(k)V(k)$, and the updated account value is therefore
$$
V(k+1) = (1 + KX(k))V(k).
$$
 With this  as the backdrop, this paper considers a new version of Kelly's problem which involves optimizing the frequency at which bets are taken. After each bet, a waiting period of~$n$ steps is enforced before updates in the wager occur. Between updates, the proceeds from the previous bets remain at risk in the spirit of ``buy and hold.''

\vspace{3mm}
When the time between bets, call it~$\Delta t$, is very small, this is viewed as the high-frequency case, and when~$\Delta t$ is large, this corresponds to ``buy and hold.'' These two extreme cases are considered in~\cite{Kuhn_Luenberger_2010} in the context of portfolio optimization with returns following a continuous geometric Brownian motion. In contrast to this work, our objective here is to analyze, in discrete time, the more general case when both the probability distribution for the returns and the time interval~$\Delta t$ between updates are arbitrary. We consider the entire range of frequencies --- from low to high.

\vspace{3mm}
In our new setting described above, the main questions which we consider are as follows: How does the optimal performance, we call it~$g_n^*$, change with the waiting period~$n \Delta t$ between updates? Does the high-frequency case~$n = 1$ always lead to the best performance? What are the effects of transaction costs? First, we provide rather complete answers to these questions for the important special case when~$X(k)$ is a Bernoulli random variable corresponding to an even-money bet with probability~$p$ of winning. This can be viewed as a simple biased coin-flipping game in the spirit of Kelly's original work and serves as an entry point for future research using a binomial lattice model for stock trading. In this part of the paper, the analysis is also extended to account for accrued interest on idle cash and transaction costs. Consistent with intuition,  when transaction costs are involved, a low frequency bettor may do  better than a high frequency bettor.

\vspace{3mm}
The latter sections of this paper focus on more general probability distributions for~$X(k)$ and two conjectures. The first conjecture is simple to state: Absent transaction costs,~$g_n^*$ is non-increasing in~$n$; i.e., $g_n^* \geq g_{n+1}^*$ for all $n$.
That is, increasing the frequency of betting can only improve performance. This raises the following questions: Assuming the conjecture is true, are the inequalities $g_{n}^* \geq g_{n+1}^*$ necessarily strict? Might it even be the case that $g_n^* = g_1^*$ for all $n$? If so, roughly speaking, this would indicate that high-frequency betting is a ``waste of time." One might as well bet at $k=0$ once without updating.

\vspace{3mm}
The second conjecture bears heavily on the questions above. To this end we work with the technical condition
$$
\mathbb{E}\left[ {\frac{1}{{1 + X(0)}}} \right] \leq 1.
$$
As explained in the sequel, this inequality is readily interpreted to be an indicator of the ``attractiveness'' of the bet. Accordingly, when this inequality is satisfied, the bet is said to be {\it sufficiently attractive}. We first prove that this condition is in fact sufficient to guarantee~$g_n^* = g_1^*$ for all~$n$. Then, we  conjecture that satisfaction of this inequality is also a necessary condition. In other words, this inequality completely characterizes the  condition under which there is no benefit associated with boosting the betting frequency. Support for this conjecture is provided by the analysis of the Bernoulli random variable  in Section~\ref{Section: MOTIVATING EXAMPLES} and numerical experiments described in Section~\ref{SECTION: TWO CONJECTURES}.

\vspace{4mm}
{\bf Kelly Criterion in a Control-Theoretic Setting}: Although our analysis here is carried out without reference to control theory, it is interesting to note the  class of problems being considered is readily reformulated in  the language of feedback systems. Such a formulation is consistent with our previous work;~e.g., see~\mbox{\cite{Hsieh_Barmisg_2015_Allerton}-\cite{Hsieh_Barmish_2017_CDC}}.
To explain further,  if we introduce the notation~$I(k)$ to denote the output of a controller which determines the ``investment'' or bet at each stage, we recognize
$
I(k) = KV(k)
$
as a linear time-invariant feedback with the gain~$K$ being the~\mbox{\it Kelly fraction}. This control-theoretic set-up is seen in~Figure~\ref{fig:Block_Diagram_KV}.  We also note that this reformulation  can be generalized to deal with the case when~$X(k)$ is a vector rather than a scalar;   see~\cite{Hsieh_Barmisg_2015_Allerton} and~\cite{Hsieh_Barmish_Gubner_2016_CDC}. For the stock market portfolio applications which we envision in our future work, this will play an important role.

%
\vspace{-1mm}
	\begin{center}
		\graphicspath{{Figs/}}
		\includegraphics[scale=0.45]{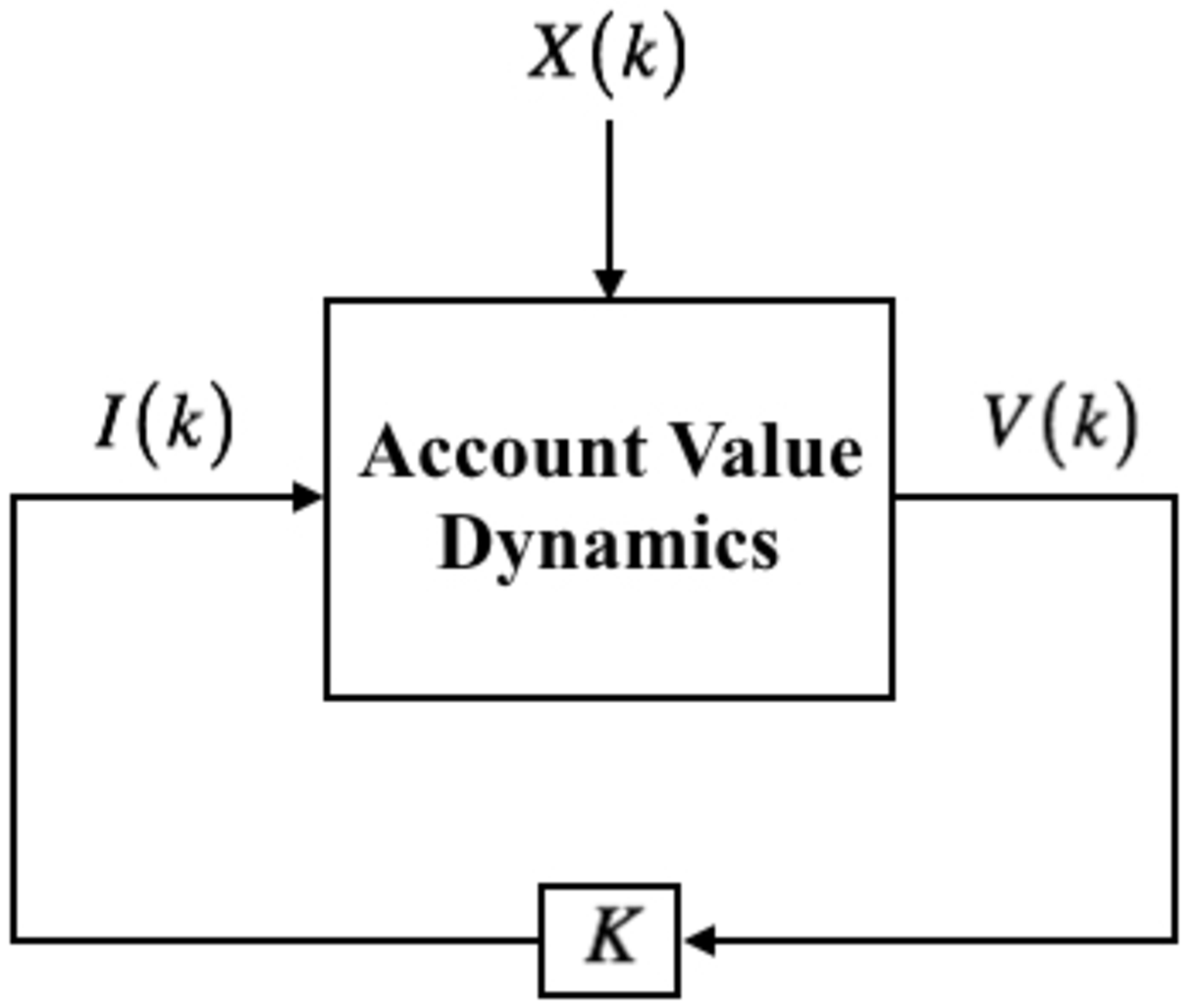}
		\figcaption{Kelly-Style Gambling Feedback Configuration}
		\label{fig:Block_Diagram_KV}
	\end{center}

\vspace{0mm}
\section{PROBLEM FORMULATION}
\label{SECTION: PROBLEM FORMULATION}
\vspace{-1mm}
Beginning with~$k = 0$, we let $X(k)$ be a sequence of i.i.d. random variables. We call~$X(k)$ the {\it returns}  and assume that
$$
X_{\min}  \leq X(k) \leq X_{\max}
$$
with $X_{\min}$ and $X_{\max}$ being support points satisfying
$$
-1 \leq X_{\min} < 0 < X_{\max} < \infty.
$$
Implicit in the inequality~$-1 \leq X_{\min}$ above is that the bettor cannot lose more than~$100\%$ of the amount bet. Now, for each integer $ n \geq 1$ denoting  the number of periods between bets, we define a new random variable
\[
\mathcal{X}_n \doteq \prod_{k=0}^{n-1} (1+X(k))-1
\]
which is induced by the~$X(k)$. This new random variable is the {\it total return} which has its own bounds obtained from~$X_{\min}$ and~$X_{\max}$. That is,
\begin{align*}
\mathcal{X}_n &\le (1+X_{\max})^n-1 \doteq \mathcal{X}_{n,\max} ; \\
\mathcal{X}_n &\ge (1+X_{\min})^n -1 \doteq \mathcal{X}_{n,\min}
\end{align*}
and it is noted that~$ \mathcal{X}_{n,\min} \geq -1$ and $ \mathcal{X}_{n,\max} < \infty $. In other words, no matter what waiting period~$n$ is used, betting losses still cannot exceed~$100\%$ and betting gains are finite. Now, for any fixed value of~$n$, we seek to maximize the expected logarithmic growth
\[
g_n(K) \doteq \frac{1}{n}\mathbb{E}[\log(1+K \mathcal{X}_n)]
\]
subject to a constraint
$
K \in {\cal K}
$
where~${\cal K}$ represents the ``betting rules" imposed by the ``house" or a broker. In this first piece of work involving frequency analysis, for simplicity we take
$
{\cal K} = [0,1]
$
in the sequel; see remarks below for further elaboration.
The associated optimal expected logarithmic growth is now obtained as
\[
g_n^* \doteq \max_{K \in \mathcal{K}}g(K)
\]
and any $K_n^* \in \mathcal{K}$ satisfying $g_n(K_n^*) = g_n^*$ is called an \mbox{\it optimal Kelly fraction}.
Given that it is trivial to show that~$\mathbb{E}[X(0) ] \leq 0$ leads to~\mbox{$K_n^* = 0$} for all $n$, that is, no  betting at all as the optimum, without loss of generality, in the sequel, our standing assumption is that
$$
\mathbb{E}[X(0)] > 0.
$$

{\bf Remarks}: The interpretation of the constraint $\mathcal{K}=[0,1]$ above is as follows: First, the requirement~$K \geq 0$ forces each bet to be non-negative. Roughly speaking this disallows taking ``opposite side'' of the bet. For example, if an even-money coin flip pays off on heads, then~$K < 0$ can be interpreted as a bet on tails which is disallowed. In the stock market, this corresponds to ``selling short.'' The second inequality associated with~${\cal K}$, namely~\mbox{$K \leq 1$}, is often called a {\it cash-financing} constraint. That is, it forces the bet size to be no more than the account value $V(k)$. In other words, this  disallows the extension of leverage from the ``house'' to the bettor. The paradigm presented in this paper remains valid for a wide variety of other choices for~${\cal K}$. Roughly speaking, extension of many of the results given here is accomplished using the notion of ``survival'' which disallows any bet that can potentially lead to~$V(k) < 0$. In the next section, although not highlighted, this survival condition is implicitly in play when transaction costs are introduced to the analysis. The reader is referred to~\cite{Hsieh_Barmish_Gubner_2016_CDC} for more details on these more general cases.

\vspace{2mm}
Another point to note is that the function~$g_n(K)$ above can be readily shown to be concave in $K$. For the case considered here with~$K$ being a scalar, the numerical maximization of~$g_n(K)$ is straightforward whether  this function is concave or not. However, for more general cases when~$K$ can be a vector of weights, concavity plays an important role because computational tractability may be an issue. The reader is referred to~\cite{Hsieh_Barmisg_2015_Allerton} and \cite{Boyd_Vandenberghe_2004} for more detailed discussion on this topic. Finally, it is important to comment on the factor~$1/n$ in the definition of expected logarithmic growth $g_n(K)$ above. It is needed so that the comparison of performance for differing values of~$n$ is ``apples to apples.'' Said another way, since~$g_n({ K})$ is associated with one bet every~$n$ periods, we adjust the expected logarithmic growth to be on a per-period basis.
If we view~$n$ as the number of bets per unit of time, then the {\it betting frequency} is given by~$1/n$.

\vspace{3mm}
\section{results for bernoulli random variables}
\label{Section: MOTIVATING EXAMPLES}
\vspace{-1mm}
Per discussion in the introduction, we
now consider the motivating case of a simple even-money  coin-flipping game. Consistent with the standing assumption that~\mbox{$\mathbb{E}[X(0)]>0$}, we assume that $p>1/2$  is the probability of heads. Now, treating the number of steps between bets as variable, the bettor contemplates the possibility of a Kelly-style adjustment of the bet size every~$n$ time steps. Between, adjustments, the money is allowed to ``ride'' with resulting profits or losses  carried over into subsequent bets and viewed as ``unrealized'' until~$n$ coin flips have been executed.

\vspace{2mm}
Indeed, since $X(k) \in \{-1, 1 \}$ is a Bernoulli random variable, our frequency-dependent method can be carried out in closed form and our results for $g_n^*$ for $n>1$ can be compared with the classical Kelly solution which, for $n=1$, is given by
$
K_1^* = 2p-1
$
with associated optimal expected logarithmic growth
\begin{align*}
	g_1^*& =  p\log(2p) + (1-p)\log(2-2p).
\end{align*}
%
The resulting  performance, as a function of frequency, is  described in the theorem below.

\vspace{3mm}
{\bf Theorem 1}: {\it For the coin-flipping game with $p>1/2$ and even-money payoff described above, the optimal Kelly fraction is
	\[
	K_n^* = \frac{2^np^n - 1}{2^n - 1}
	\]
and the associated optimal expected logarithmic growth is given by
\[
g_n^* = p^n\log{p} + \left(\frac{1-p^n}{n}\right)\log\left[{\frac{1-p^n}{2^n-1}}\right]+ \log{2}.
\]}%
%
{\bf Proof}:
We first observe that the total return is given by
$$
{\cal X}_n = 2^n-1\;\;\text{with probability}\;\; p^n
$$
and
$$
{\cal X}_n = -1\;\;\text{with probability}\;\; 1-p^n.
$$
We now calculate the expected logarithmic growth
{\small
\begin{eqnarray*}
	g_n({ K}) & \doteq& \frac{1}{n} \mathbb{E} \left[\log(1 + { K}{\cal X}_n)\right] \\ [5pt]
	& = & \frac{1}{n}\left[p^n\log(1 + (2^n-1){ K}) + (1-p^n)\log(1-{ K})\right]
\end{eqnarray*}
}
and find its derivative
$$
\frac{\partial{g_n}}{\partial{ K}} = \frac{1}{n}\left[\frac{p^n(2^n - 1)}{1 + (2^n-1){ K}} -\frac{1-p^n}{1-{ K}}\right].
$$
Setting this to zero, we obtain a unique candidate for the maximizer, call it~${ K} = { K}_n^*$, given by
$$
{ K}_n^* = \frac{2^np^n - 1}{2^n - 1}
$$
which is feasible; i.e.,~${ K}_n^*\in [0,1]$.
Now using the fact that~$g_n(K)$ is concave, in combination with uniqueness of the zero-derivative point above, it follows that $K_n^*$ is the global maximizer; e.g., see \cite{Boyd_Vandenberghe_2004}.
Finally, to complete our analysis, we substitute~${ K}_n^*$ into~$g_n({ K})$. A  lengthy but straightforward calculation leads to
$$
g_n^* = p^n\log{p} + \left(\frac{1-p^n}{n}\right)\log\left[{\frac{1-p^n}{2^n-1}}\right]+ \log{2}. \;\;\;\square
$$
\vspace{-9mm}
\begin{center}
	\graphicspath{{figs/}}
	\includegraphics[scale=0.4]{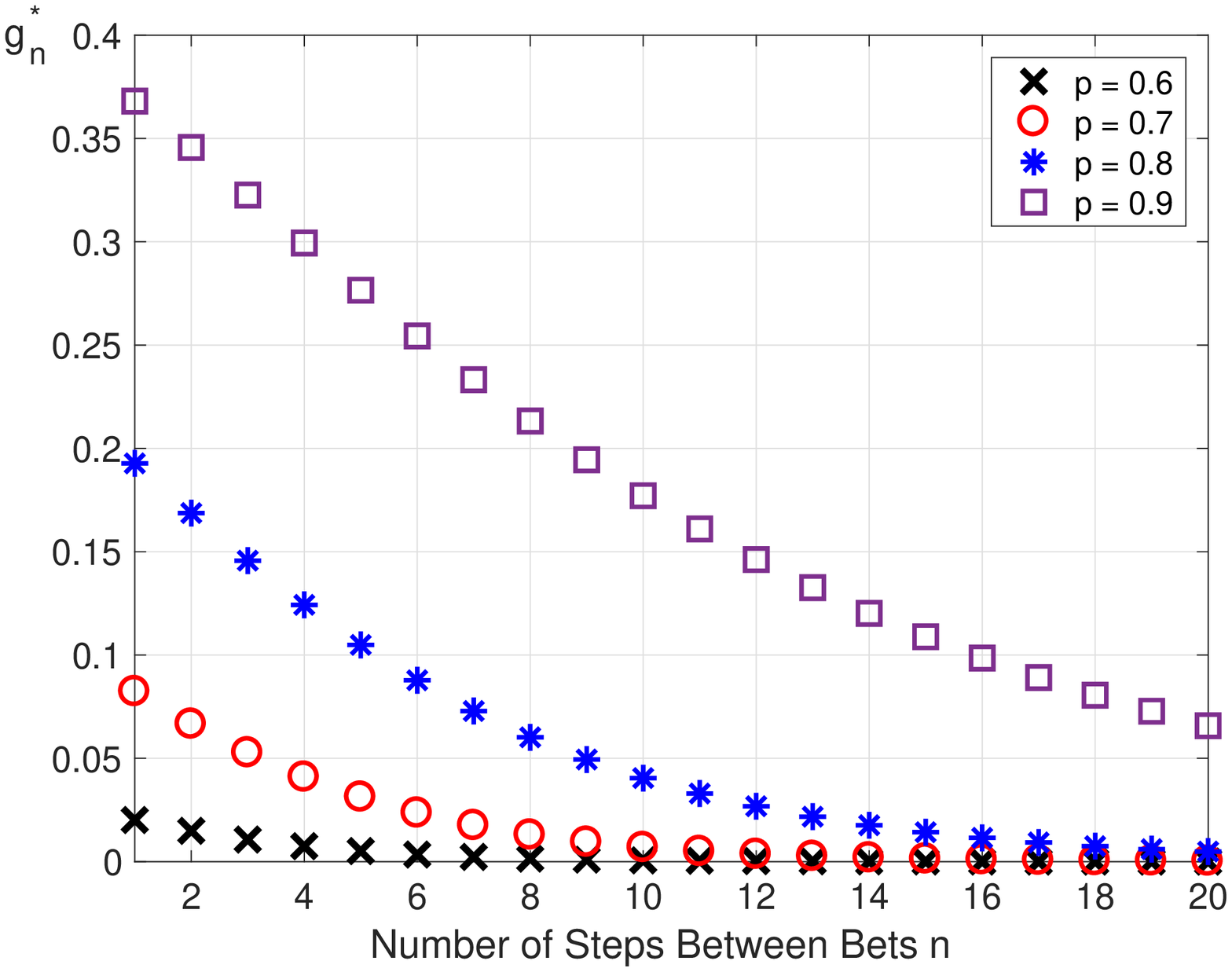}
	\figcaption{Optimal Expected Logarithmic Growth Versus $n$ }
	\label{fig:Optimal Expected Logarithmic Growth Versus n}
\end{center}

\vspace{3mm}
{\bf Remarks}:
To more clearly demonstrate frequency dependence of the performance resulting from the formulae above, in Figure~\ref{fig:Optimal Expected Logarithmic Growth Versus n}, for various $p$ values, we  plot  $g_n^*$ as a function of the waiting period $n$. Noting that~$g_n^*$ is  non-increasing in~$n$ which suggests that betting ``faster is better." In the next section, we provide conjectures regarding the extent to which such results are generalizable to more general probability distributions for the $X(k).$



\vspace{3mm}
\section{transaction costs and accrued interest}
\label{SECTION: accrued interest and transaction costs} 
\vspace{-1mm}
We now indicate how the problem formulation in this paper is readily extended to include consideration of transaction costs and accrued interest on idle cash.
As far as transaction costs are concerned, we consider the case when a percentage charge~$0 <  \varepsilon \leq  1$ is assessed by the ``house'' each time the bet is updated. If this occurs at stage~$k$, the charge is~$\varepsilon K V(k)$. Common sense tells us that as~$\varepsilon$ increases, the benefits of high-frequency betting are negated. This is confirmed in the example to follow.
 In addition, for the case of idle cash, let~$r$ denote the per-period rate of interest. With~$n$ steps between updates of the bet, the accrued interest from~$k = 0$ to~$k = n$ is~$(1+r)^n(1-K - \varepsilon K)V(0)$ and it is natural to anticipate that the optimal fraction~$K_n^*$ should decrease as~$r$ increases; that is, the incentive to tie up money in a risky gamble is reduced when a riskless alternative is available. 

\vspace{3mm}
For the sake of brevity, we simply summarize how the two considerations above are handled in an extended formulation. Suffice it say, the effects of interest and transaction costs can be ``lumped into'' the returns. After a straightforward calculation, we obtain modified expected logarithmic growth given by
$$
g_n(K) = \log(1+r) + \frac{1}{n} \mathbb{E}[\log(1 + K\tilde{{\cal X}}_n)]
$$
with 
$$
\tilde{{\cal X}}_n \doteq \frac{1 + {\cal X}_n  }{(1+r)^n}-\varepsilon   - 1.
$$
In summary, in the high level formulation, whenever convenient, there is no loss of generality taking~$\varepsilon = 0$ and~$r = 0$.

%

\vspace{3mm}
{\bf Example: Transaction Cost Considerations}: To illustrate how transaction costs negate the performance effects of high-frequency betting, we revisit the even-money Bernoulli case  considered in the previous section. We now demonstrate that~$n=1$ is no longer optimal as the transaction costs grow. To this end, per discussion above, we now take~$0 < \varepsilon <1 $ to be a  percentage  transaction cost which is incurred by the bettor each time the bet size is updated. 
Note that when these costs costs are included, the possibility exists that for~$K> 1/(1+\varepsilon)$, the situation $V(k)<0$ will arise. To avoid such ambiguities, in the analysis to follow, we work with the subset of $\mathcal{K}$ where the ``survival" issue on~$K$ is involved; i.e., the fraction
\[
K \in \mathcal{K}_\varepsilon \doteq \left[ {0,\;\frac{1}{{1 + \varepsilon }}} \right].
\]
Now, to isolate these transaction cost effects, we take interest rate $r=0$ and find $K_n^* \in \mathcal{ K}_\varepsilon$ maximizing the expected logarithmic growth given by
\begin{align*}
g_{n,\varepsilon}(K) &= \frac{1}{n}\mathbb{E}\left[ \; {\log \left( {1   + K( \mathcal{X}_n - \varepsilon)} \; \right)} \right].
\end{align*}
Note that $g_{n,\varepsilon}(K)$ is concave in $K$.

\vspace{3mm}
{\bf Theorem 2}: {\it Consider the coin-flipping game described above with  $p>1/2$, even-money payoff,  and  transaction cost~\mbox{$0 < \varepsilon <1$}. Then, with
\[
p \leq {p_\varepsilon } \doteq {\left( {\frac{{1 + \varepsilon }}{{{2^n}}}} \right)^{1/n}},
\]
the optimal Kelly fraction is $K_{n,\varepsilon}^*=0$ with associated expected logarithmic growth $g_{n,\varepsilon}^*=0$. Otherwise, the optimum is given by
\[K_{n,\varepsilon }^* = \frac{{{2^n}{p^n} -\varepsilon - 1  }}{{(\varepsilon  + 1)\left( {{2^n} - \varepsilon  - 1} \right)}}\]
with the associated expected logarithmic growth given by
\begin{align*}
g_{n,\varepsilon }^* = {p^n}\log \left( {\frac{{{2^n}{p^n}}}{{\varepsilon  + 1}}} \right) + \left( {1 - {p^n}} \right)\log \left( {\frac{{{2^n}\left( {1 - {p^n}} \right)}}{{{2^n} - \varepsilon  - 1}}} \right).
\end{align*}

}

\vspace{2mm}
{\bf Proof}:
To begin, note that the expected logarithmic growth for~${K} >0$ is given by
\begin{align*}
g_{n,\varepsilon }({K}) &= p^n \log \left( {1 - \varepsilon {K} + {K}\left( {{2^n} - 1} \right)} \right)
\\
& \hspace{5mm} + \left( {1 - p^n} \right)\log \left( {1 - \varepsilon {K} - {K}} \right).
\end{align*}
Now taking the first derivative and setting it equal to zero, we obtain a candidate for the  optimal Kelly fraction, call it~\mbox{$K=K_{n,\varepsilon}^*$}  given above.
To complete the proof, we first show that $K_{n,\varepsilon}^* \in \mathcal{ K}_\varepsilon$. Indeed,
observing that $2^n - \varepsilon -1 >0$ and $1/2 < p \leq 1 $, it is straightforward to see that
\[
 K_{n,\varepsilon}^* \le \frac{1}{{1 + \varepsilon }}.
\]
In addition, for $p \geq p_{\varepsilon}$, we have $K_{n, \varepsilon}^* \geq 0.$
Now using the fact that the zero-derivative point~${K}_{n,\varepsilon}^*$ is unique, in combination with the fact that $g_{n,\varepsilon}(K)$, arguing as in the proof of Theorem 1, is concave, it follows that~$K_{n,\varepsilon}^*$ is the global maximizer. Therefore, substituting~$K_{n,\varepsilon}^*$ into $g_{n, \varepsilon}(K)$, a straightforward calculation leads
 to~$g_{n,\varepsilon}^*$ as required. \hspace{5mm} $\square$

%
%

\vspace{-3mm}
\begin{center}
	\graphicspath{{figs/}}
	\includegraphics[scale=0.40]{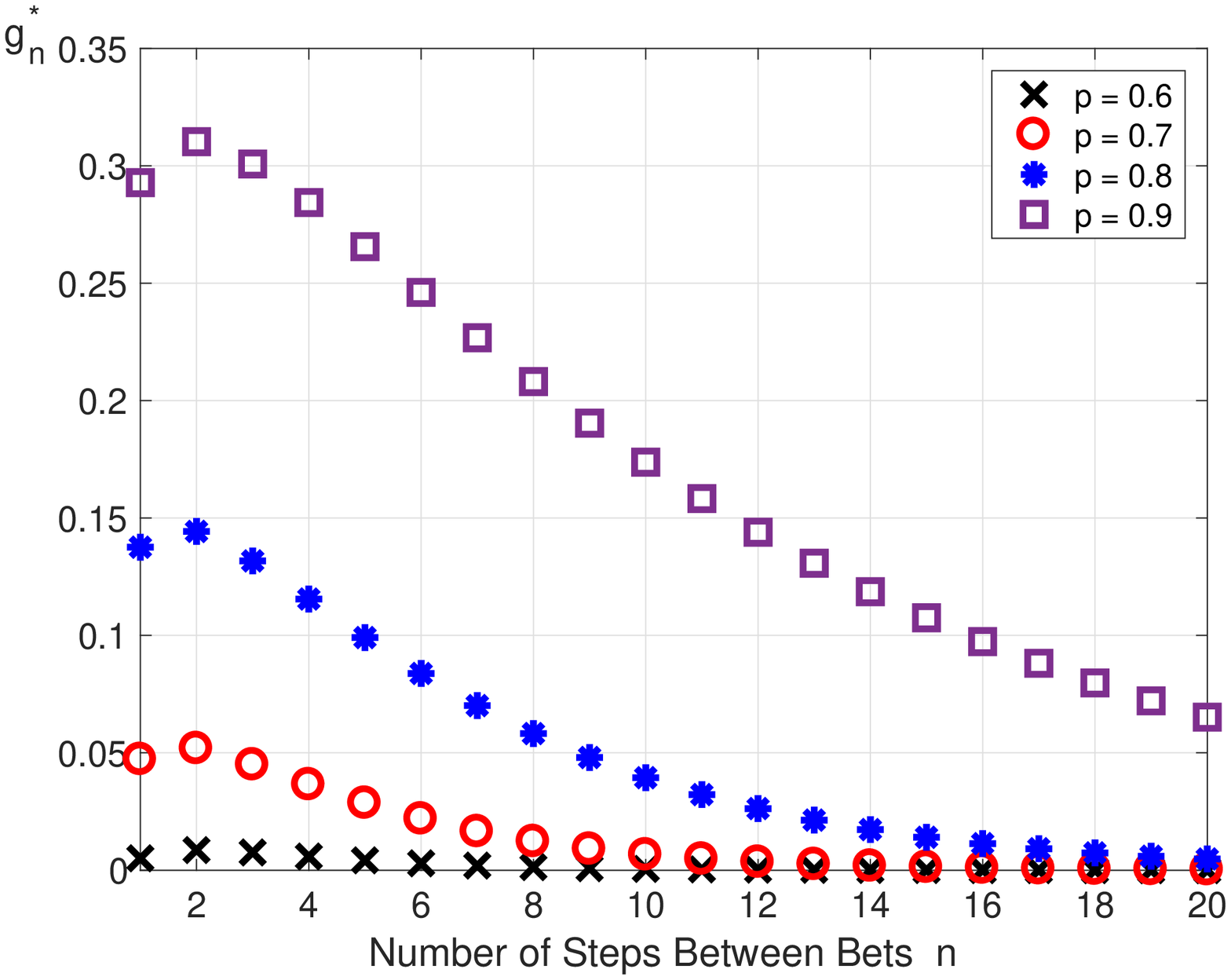}
	\figcaption{Transaction Cost Effects on $g_n^*$ }
	\label{fig.Transaction Cost Effects in g_n^*}
\end{center}

\vspace{2mm}
{\bf Remarks}: Note that when $\varepsilon \to 0$, we obtain the original zero-transaction cost result in the previous theorem as a special case; i.e.,
\begin{eqnarray*}
  \mathop {\lim }\limits_{\varepsilon  \to 0} K_{n,\varepsilon }^* &=& \mathop {\lim }\limits_{\varepsilon  \to 0} \frac{{{2^n}{p^n} - 1 - \varepsilon }}{{(\varepsilon  + 1)\left( {{2^n} - \varepsilon  - 1} \right)}} \hfill \\[5pt]
   &=& \frac{{{2^n}{p^n} - 1}}{{{2^n} - 1}} \\[5pt]
   &=& K_n^*.
\end{eqnarray*}

To illustrate use of the results above, for  $p=0.6,0.7,0.8$ and~$0.9$  and~$\varepsilon = 0.1$, we plot the optimal expected logarithmic growth $g_n^*$ versus~$n$ in Figure~\ref{fig.Transaction Cost Effects in g_n^*}.
From that figure, we see that the optimal waiting period is $n^*=2.$ In other words, the optimum is achieved by waiting two time units rather than betting at each step.

\vspace{3mm}
 \section{two conjectures and supporting evidence}
 \label{SECTION: TWO CONJECTURES}
 \vspace{-1mm}
 In this section, we provide two conjectures which are suggested by the  results in the previous sections and numerical experiments.
Our first conjecture, roughly speaking,  says that absent transaction costs,  betting with higher frequency leads to better performance.   Roughly speaking this says that ``faster is better." Indeed, this is clearly seen to be true for the special case involving a simple coin-flipping example. The conjecture below pertains to the general case formulated in  Section~\ref{SECTION: PROBLEM FORMULATION}.


 \vspace{3mm}
{\bf Conjecture 1}: {\it The optimal expected logarithmic growth $g_n^*$ is non-increasing in $n$; i.e., for all~\mbox{$n \geq 1,$}
$
	g_{n}^* \ge g_{n+1}^*.
$	
}

\vspace{3mm}
{\bf Supporting Evidence For Conjecture 1}: The example which we now analyze is a generalization of the Bernoulli random variable scenario in Section~\ref{Section: MOTIVATING EXAMPLES}. Indeed, we now consider a coin-flipping gamble described by
$$
X(k) \in \{-\gamma, \gamma\}
$$
with~\mbox{$0 < \gamma \leq 1$} and $p > 1/2$. In comparison with our earlier scenario with~$\gamma = 1$,  the analysis becomes more complicated for the following reason: With~$0< \gamma < 1$, the total return~$\mathcal{X}_n$ no longer has a two-point probability mass function. In this  case, for~$i = 0,1,2,\ldots,n$ we obtain point masses located at
$$
x_i = (1+\gamma)^i(1-\gamma)^{n-i} -1
$$
with associated probabilities
\[
{p_i} = P\left( \mathcal{X}_n = {x_i} \right) = \left( \begin{gathered}
  n \hfill \\
  i \hfill \\
\end{gathered}  \right){p^i}{\left( {1 - p} \right)^{n - i}}.
\]
This leads to the formula
\[
{g_n}\left( K \right) = \frac{1}{n}\sum\limits_{i = 0}^n {{p_i}\log \left( {1 + K{x_i}} \right)}
 \]
which we maximized numerically with respect to~\mbox{$K \in {\cal K} = [0,1]$.} To study the conjecture above, we carried out the optimization for many values of~$\gamma$. The case~$\gamma = 1/2$, which we now describe, typifies our findings.
To  clearly see the difference between the~$g_1^*$ and~$g_n^*$, we plotted $e_n^* \doteq g_1^* -g_n^*$ as a function of $p$  in Figure~\ref{fig.Error Plot_for_UnevenPayoff}.
Consistent with Conjecture~1, we see~$g_n^* >  g_{n+1}^*$ for~$p \leq 0.75$.

\vspace{3mm}
More interestingly, however, is our observation that for~\mbox{$p \ge 0.75$}, the inequality relating $g_n^*$ to $g_{n+1}^*$ ceases to be strict. That is, once the bet becomes ``sufficiently attractive,'' there does not appear to be any benefit associated with high-frequency betting. For example, the expected logarithmic growth betting every ten steps is the same as what one obtains betting every step. This observation paves the way for our second conjecture which bears on the following roughly-stated question: Under
what conditions is increasing the betting frequency a waste of time? To address this issue   a technical condition,  which we call the  \mbox{\it sufficient attractiveness} inequality, is now introduced.

\vspace{3mm}
	\begin{center}
		\graphicspath{{figs/}}
		\includegraphics[scale=0.85]{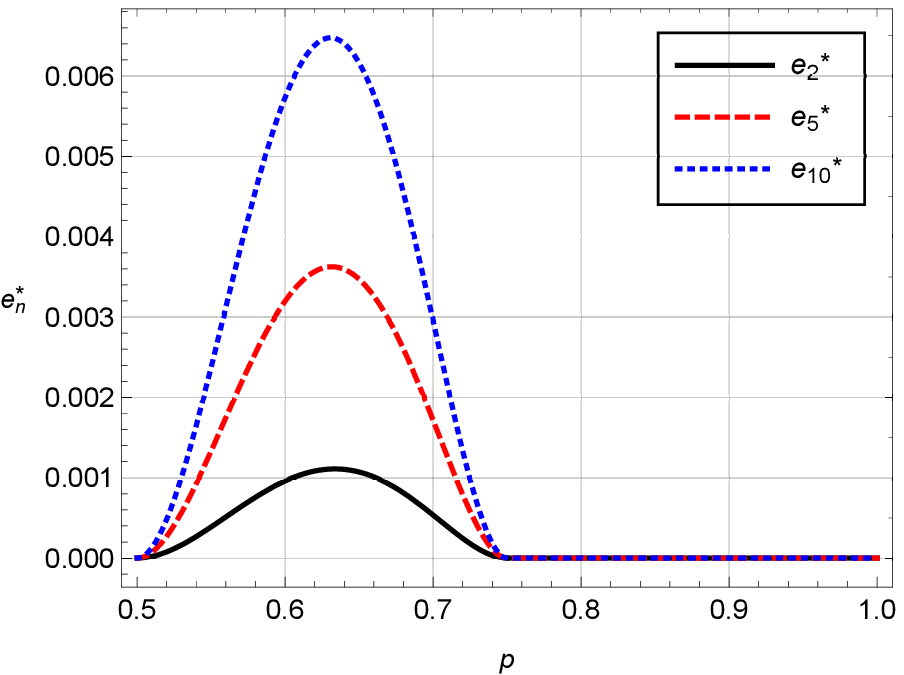}
		\figcaption{Plot of Error $e_n^*$ Versus $p$}
		\label{fig.Error Plot_for_UnevenPayoff}
	\end{center}

\vspace{3mm}
{\bf Definition}: The i.i.d. random variables $X(k)$ are said to satisfy the {\it sufficient attractiveness} inequality if
\[
\mathbb{E} \left[ {\frac{1}{{1 + X\left( 0 \right)}}} \right] \le 1.
\]

\vspace{2mm}
{\bf Remarks}:  The use of the descriptor ``sufficient attractiveness'' above is consistent with our experience with a wide variety of examples which indicates the following:  This inequality tends to be satisfied when the gamble becomes highly favorable to the bettor. As a first illustration, for the extreme case with~$X(0) > 0$ with probability one, we have an arbitrage, and the sufficient attractiveness inequality  above is trivially satisfied. Indeed, we now argue that the smallness of $\mathbb{E}[1/(1+X(0))]$ is desirable. To establish this, put
$$
\theta \doteq \mathbb{E}\left[ {\frac{1}{{1 + X(0)}}} \right].
$$
Then, since $f(z) = 1/(1+z)$ is convex, applying the Jensen's inequality leads to
\[
\frac{1}{{1 + \mathbb{E}\left[ {X(0)} \right]}} \leq \mathbb{E}\left[ {\frac{1}{{1 + X(0)}}} \right],
\]
which, for $0 < \theta \leq 1$ implies that
\[
\mathbb{E}[X(0)] \geq \frac{1}{\theta} -1.
\]
This inequality tells us that as  $\theta$ decreases to zero,  $\mathbb{E}[X(0)]$ increases without bound. In other words,  smallness of $\mathbb{E}[1/(1+X(0))]$ is desirable. The inequality also tell us that sufficient attractiveness implies $\mathbb{E}[X(0)]\geq 0.$

\vspace{2mm}
More realistically, for the case considered in the previous section with~$X(k) = \gamma$ with probability~$p$ and  $X(k) = -\gamma$ with probability~$1-p$, it is easy to verify that the sufficient attractiveness inequality is satisfied if and only if
\[
p \ge \frac{{1 + \gamma }}{2}.
\]
For example, for the even-money payoff with~$\gamma = 1/2$, we require~$p  \geq 0.75$. Note that for the extreme case~\mbox{$\gamma=1$}, we require $p=1$ and all $g_n^*$ are equal. More generally, for Bernoulli problems with \mbox{$X(k) \in \{X_{\min},X_{\max}\}$} and~$P(X(k)=X_{\max})=p$, the sufficient attractiveness inequality reduces to
\[
p \ge  \frac{{  {|X_{\min }|}\left( { 1+ {X_{\max }} } \right)}}{{{X_{\max }} - {X_{\min }}}}.
\]
This  tells us that the winning probability $p$ should be higher than some threshold in order for the bet to be recognized as sufficiently attractive. Below we provide a  result which serves as a feeder for our next conjecture.

\vspace{3mm}
{\bf Sufficient Attractiveness for Uniform Distribution}: 
We now provide more detail for the case of the uniform distribution. 
This example is interesting because the sufficient attractiveness inequality can be analyzed in closed form. Suppose the i.i.d. returns~$X(k)$ are governed by the uniform distribution on~$[a,b]$ with~\mbox{$-1<a<0$} and $b>0$. To first study  the sufficient attractiveness inequality, we calculate
\begin{align*}
\mathbb{E}\left[ {\frac{1}{{1 + X(0) }}} \right]
&= \frac{1}{{b - a}}\log \left( {\frac{{1 + b}}{{1 + a}}} \right).
\end{align*}
Now, for $a$ in its allowed range, we define
\begin{align*}
{b_{\min }}(a) &= \min \left\{ {b > 0: \mathbb{E}\left[ {\frac{1}{{1 + X\left( 0 \right)}}} \right] \leq 1} \right\}.
\end{align*}
Note that the sufficient attractiveness inequality for the uniform distribution case is easily seen to be equivalent to
$$	
\frac{e^a}{1+a} \leq \frac{e^b}{ 1+b}.
$$	
If we put $f(t)=e^t/(1+t)$, then the inequality above says~\mbox{$f(a) \leq f(b)$}. Since~$-1 < a < 0 < b$ is assumed  in the uniform case, the bet is sufficiently attractive if and only if~\mbox{$b \geq   {b_{\min }}(a)$}, where~$  {b_{\min }}(a) > 0$ solves~\mbox{$f(b)=f(a)$}, which is easily accomplished in Matlab. A formula for~$b_{\min}(a)$ is provided  in the remark below.

\vspace{3mm}
{\bf Remarks}: The description of $b_{\min}(a)$ above can also be given in terms of the Lambert function. Namely,
\[
b_{\min}(a) = -1 - W_{_{-1}} \left( { - \left( {1 + a} \right){e^{ - \left( {1 + a} \right)}}} \right)
\]
where  $W_{_{-1}}(\cdot)$ is the Lambert function with lower branch satisfying $W_{_{-1}}(\cdot) \leq -1$.
Figure~\ref{fig:SAbet_uniform} provides a plot for  $b_{\min} (a)$ versus $|a|$ for~\mbox{$a \in (-1,0)$}. Note that as $a$ decreases, the figure tells us what  value of $b$ is needed to sustain  sufficient attractiveness.

\vspace{-1mm}
\begin{center}
	\graphicspath{{figs/}}
	\includegraphics[scale=0.43]{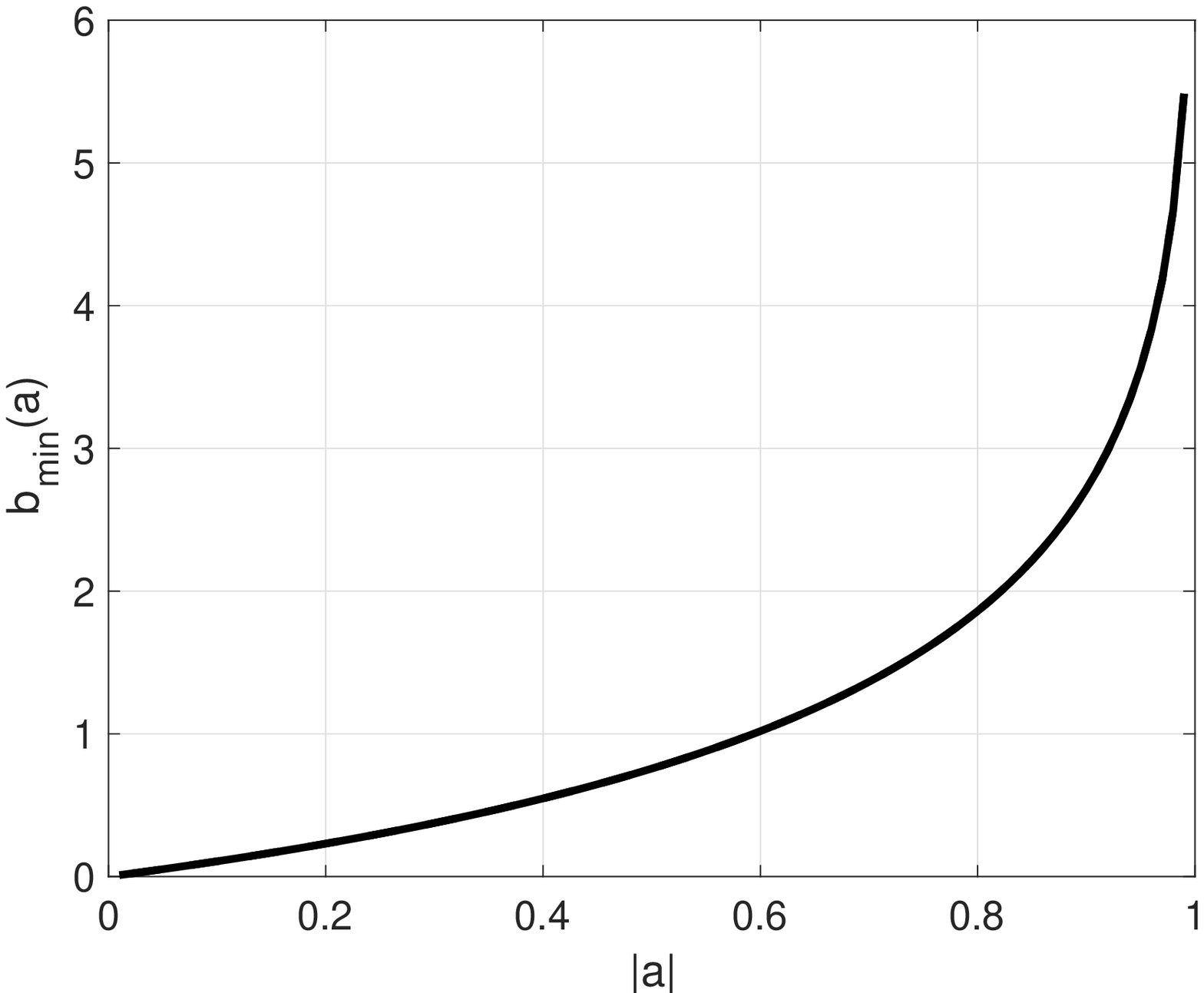}
	\figcaption{Sufficient Attractiveness: $b_{\min}(a)$ Versus $|a|$}
	\label{fig:SAbet_uniform}
\end{center}

\vspace{3mm}
{\bf Sufficiency Theorem}: {\it  Satisfaction of the sufficient attractiveness inequality guarantees that  the condition
$
	g_1^* = g_n^*
$
holds for all~$n\ge 1$.
}

\vspace{3mm}
{\bf Proof}: Assuming that the sufficient attractiveness inequality holds, we first claim that $g_n(K)$ is nondecreasing. Beginning with
\begin{align*}
\frac{d}{dK}g_n(K)
&= \frac{d}{{dK}} \mathbb{E}\left[ {\log (1 + K{ \mathcal{X}_n})} \right]
\end{align*}
and noting that $\mathcal{X}_n$ is bounded, results in probability theory, for example, see~\cite{Billingsley_1995} and~\cite{Gubner_2006}, allow us to commute the differentiation and expectation operators above. Hence,
%
\begin{eqnarray*}
\frac{d}{{dK}}\mathbb{E}\left[ {\log (1 + K{\mathcal{X}_n})} \right]
&=& \mathbb{E}\left[ {\frac{d}{{dK}}\log (1 + K{\mathcal{X}_n})} \right]\\[5pt]
&=& \mathbb{E}\left[ {\frac{{{\mathcal{X}_n}}}{{1 + K{\mathcal{X}_n}}}} \right].
\end{eqnarray*}
Now noting that the inequality
\begin{align*}
\frac{z}{1+Kz} &\geq 1-\frac{1}{1+z}
\end{align*}
holds for all~$K \in [0,1]$ and all $z > -1$, we obtain
\begin{eqnarray*}
\frac{d}{{dK}}\mathbb{E}\left[ {\log (1 + K{\mathcal{X}_n})} \right]
&\geq & 1 - \mathbb{E}\left[ {\frac{1}{{1 + {\mathcal{X}_n}}}} \right].
\end{eqnarray*}
Using the fact that the $X(k)$ are i.i.d., we observe that
\begin{eqnarray*}
\mathbb{E}\left[ {\frac{1}{{1 + {\mathcal{X}_n}}}} \right]
& = & \mathbb{E}\left[ {\prod\limits_{k = 0}^{n - 1} {\frac{1}{{1 + X\left( k \right)}}} } \right]\\[5pt]
&=& {\left( {\mathbb{E}\left[ {\frac{1}{{1 + X\left( 0 \right)}}} \right]} \right)^n}\\[5pt]
& \leq& 1.
\end{eqnarray*}
Now combining this with the previous inequality, we have
\[
\frac{d}{dK} g_n(K) \ge 0
\]
which shows that $g_n(K)$ is non-decreasing in $K$. Hence~$g_n(K)$ is maximized at $K=1.$ Therefore,
for all~$n$, we have $g_n^* = g_n\left( 1 \right)$.
It only remains to observe that
\begin{eqnarray*}
g_n\left( 1 \right) &=& \frac{1}{n} \mathbb{E}[\log \left( {1 + { \mathcal{X}_n}} \right)]\\ [5pt]
&=& \frac{1}{n}\sum\limits_{k = 0}^{n - 1} {\mathbb{E} [\log \left( {1 + X\left( k \right)} \right)}] \\[5pt]
&=& \mathbb{E} [\log \left( {1 + X\left( 0 \right)} \right)] \\[5pt]
& =& g_1^*
\end{eqnarray*}
which completes the proof. \;\;\;  $\square$

\vspace{3mm}
{\bf Remarks}: The Sufficiency Theorem tells us that when  a bet is attractive enough, the buy-and-hold bettor can achieve the same betting performance as the high-frequency bettor.
It is also worth mentioning that the degenerate case $\mathbb{E}[X(0)]=0$, not covered by our initial assumptions, is easily shown to lead to all $g_n^*=0$ too; i.e., the optimum is $K_n^* = 0$ which corresponds to no betting.

\vspace{3mm}
{\bf Conjecture 2}: {\it Satisfaction of the sufficient attractiveness inequality is necessary for the condition~$
	g_1^* = g_n^*
	$
	to hold for all $n\ge 1$.
}

\vspace{3mm}
{\bf  Supporting Evidence for Conjecture 2}:
It is  first noted that we carried out a number of preliminary numerical experiments with results which support the conjecture for commonly used distributions such as Bernoulli, uniform, triangular, beta and truncated normal. 
That is, with the sufficient attractiveness inequality satisfied, for many cases, we carried out the required optimizations to find the $g_n^*$. To numerical accuracy we found $g_n^* \approx g_1^*$
 for the ranges of $n$ which we considered.To provide further supporting evidence for the conjectured necessity of sufficient attractiveness, observe in Figure~\ref{fig.Error Plot_for_UnevenPayoff} that for~$p \geq 0.75$, we have~$g_{n}^* =  g_1^*$ for~$n=2,5$, and $10$. For the  uniform distribution case, although not shown in the figure,  our numerical experiments also show that $g_n^* = g_1^*$ for~\mbox{$n \leq 20$}.

\vspace{4mm}
\section{CONCLUSION AND FUTURE WORK}
\label{CONCLUSION AND FUTURE WORK}
\vspace{-1mm}
In this paper, we studied the problem of optimizing the betting frequency in a dynamic setting using Kelly's celebrated expected logarithmic growth criterion as the performance metric. We provided a detailed analysis for the important special case when~$X(k)$ is a Bernoulli random variable with even-money payoff. We also extended our analysis to include accrued interest on idle cash and transaction costs. The results for the Bernoulli case support our first conjecture that the optimal expected logarithmic growth $g_n^*$ is a non-increasing in~$n$.

\vspace{2mm}
We also investigated conditions  under which betting with arbitrarily low  frequency can still achieve the same  performance as betting with very high-frequency. To this end, we showed that satisfaction of the  sufficient attractiveness inequality assures that~$g_n^* = g_1^*$ for all~$n$. Subsequently, we  conjectured that the sufficient attractiveness is also necessary and provided numerical evidence in support.

\vspace{2mm}
Regarding further research, one obvious continuation would be to study the two conjectures in greater detail.  Another  problem for future research involves extension of the results in this paper to a portfolio scenario with  many correlated random variables; i.e., we take $X(k)$ to be a vector rather than the scalar considered here.
For this more general case with~$K$ being a vector of weights, we envision concave programming playing an important role since computational tractability can become a significant issue when multi-dimensional optimization is performed.

\end{document}